\documentclass{amsart}
\usepackage{amssymb}
\usepackage{latexsym}
\usepackage{amsmath}
\usepackage{euscript}
\usepackage{graphics}
\usepackage[all]{xy}

\usepackage[leqno]{amsmath}

\usepackage{bbm}

\usepackage{dsfont}

\usepackage{graphicx}

      \def\dC{{\mathbb C}}

      \def\dR{{\mathbb R}}

\newcommand{\be}{\begin{equation}}
\newcommand{\ee}{\end{equation}}
\newcommand{\ba}{\begin{eqnarray}}
\newcommand{\ea}{\end{eqnarray}}
\newcommand{\baa}{\begin{eqnarray*}}
\newcommand{\eaa}{\end{eqnarray*}}
\newcommand{\bb}{\begin{thebibliography}}
\newcommand{\eb}{\end{thebibliography}}

\newcommand{\bi}[1]{\bibitem{#1}}
\newcommand{\lab}[1]{\label{#1}}
\newcommand{\re}[1]{(\ref{#1})}

\newcounter{my}
\newcommand{\he}%
   {\stepcounter{equation}\setcounter{my}%
   {\value{equation}}\setcounter{equation}0%
   }%
\newcommand{\she}%
   {\setcounter{equation}{\value{my}}%
    }%

\renewcommand\t{\tilde}
\def\sign{\operatorname{sign}}

\newtheorem{theorem}{Theorem}[section]

\newtheorem{lemma}[theorem]{Lemma}

\theoremstyle{definition}

\newtheorem{remark}[theorem]{Remark}

\numberwithin{equation}{section}

\begin{document}

\title[CMV and -1 Jacobi polynomials]
{CMV matrices and little and big -1 Jacobi polynomials}

\author{Maxim Derevyagin}
\author{Luc Vinet}
\author{Alexei Zhedanov}

\address{Department of Mathematics MA 4-5, Technische Universit\"at Berlin, Strasse des 17. Juni 136, D-10623 Berlin, Germany}

\address{Centre de recherches math\'ematiques,
Universit\'e de Montr\'eal, P.O. Box 6128, Centre-ville Station,
Montr\'eal (Qu\'ebec), H3C 3J7,Canada}

\address{Donetsk Institute for Physics and Technology\\
83114 Donetsk, Ukraine \\}

\begin{abstract}
We introduce a new map from polynomials orthogonal on the unit
circle to polynomials orthogonal on the real axis.
This map is closely related with the theory of CMV matrices. It
contains an arbitrary parameter $\lambda$ which leads to a linear
operator pencil. We show that the little and big -1 Jacobi
polynomials are naturally obtained under this map from the Jacobi
polynomials on the unit circle.
\end{abstract}

\keywords{Classical orthogonal polynomials, Jacobi polynomials,
little and big -1-Jacobi polynomials. AMS classification: 33C45,
33C47, 42C05}

\maketitle

\section{Introduction}
In this paper we propose a new map from polynomials orthogonal on
the unit circle (OPUC) to polynomial orthogonal on the real line.
This map can be considered as a one-parameter generalization of the
Delsarte-Genin \cite{DG1} map from OPUC to symmetric orthogonal
polynomials on an interval. We derive this map from
elementary properties of the CMV matrices. The CMV matrices are
known to play a crucial role in the theory of OPUC \cite{Simon},
\cite{Watkins}, \cite{CMV}.

This new map contains an additional arbitrary parameter $\lambda$
which leads to infinitely many families of polynomials
$P_n(x;\lambda)$ orthogonal on the real line. In general, for
arbitrary real values of $\lambda$ these polynomials are
orthogonal on a union of two distinct intervals of the real line.

Starting from given explicit sets of OPUC, one thus constructs
families of polynomials orthogonal on the real line
$P_n(x;\lambda)$. The main result of this paper consists in the
following. Starting from the Jacobi polynomials on the unit circle, we
construct a one-parameter family of polynomials $P_n(x;\lambda)$
orthogonal on the union of two intervals of the real line.

We identify these polynomials with the big -1 Jacobi polynomials
discovered recently in \cite{VZ_big}. For a particular value
$\lambda=1$ of the parameter we obtain the little -1 Jacobi
polynomials found in \cite{VZ_little}.

The underlying spectral problem has the following form
\begin{equation}\label{LinPen2par}
(L+\lambda M-xI)\vec q(x,z)=0,
\end{equation}
where $L$ and $M$ are tridiagonal matrices of a special form
($L$ and $M$ have zeros on the super- and sub-diagonals).

If we fix $\lambda$, the problem~\eqref{LinPen2par} is a classical eigenvalue problem
for the tridiagonal matrix $K(\lambda)=L+\lambda M$. Thus, it gives rise to a system of orthogonal polynomials,
Pad\'e approximants and moment problems.

If we fix $x$,  the problem~\eqref{LinPen2par} is a generalized eigenvalue problem
for the two matrices $M$ and $L-xI$. It generates a system of functions rational in $x$ that
are related to the theory of bi-orthogonal rational functions~\cite{ZheBRF}.
Moreover, it has recently been shown that such generalized eigenvalue problems
appear in rational interpolation~\cite{BDZh},~\cite{DZh},
Nevanlinna-Pick problems~\cite{Der10},~\cite{DZh},
integrable systems~\cite{SZ},~\cite{VSh}, and mechanics~\cite[Section~2.4]{Gladwell}
(see also~\cite{DG2,DG3}).

A priori, the general theory  associated to the problem~\eqref{LinPen2par} is rather involved.
That is why we make some assumptions on the matrices $L$ and $M$. Namely, we assume that $L^2=I$ and $M^2=I$.
 Thus, the main idea is to identify systems of orthogonal polynomials that have a hidden
biorthogonality associated to a parameter.

Note that the special case $x=0$ leads to the  Schur linear pencil
$L + \lambda M$ considered by Watkins \cite{Watkins}.

Finally it should be mentioned that particular cases of~\eqref{LinPen2par}
surprisingly arise in the theory of martingale
polynomials~\cite{BMW},~\cite{ST98}, two-state free Brownian
motions~\cite{Ans10}, and Gaussian processes~\cite{Ker98}. Let us also signal references~\cite{ER}
and~\cite{Stanton} where similar
pencils appear in the investigation of exactly solvable kinetic models
and the corresponding combinatorial problems.

The paper is organized as follows.

In Section 2, we recall elementary  facts concerning OPUC and CMV
matrices.

In Section 3, we introduce the Jacobi (3-diagonal) matrix $J$
starting from two unitary elementary matrices $L$ and $M$. A new
family $Q_n(x)$ of orthogonal polynomials on the real line is
associated with the matrix $J$.

In Section 4, basic results on the Delsarte-Genin map from
symmetric orthogonal polynomials $S_n(x)$ to OPUC are presented.

In Section 5, we introduce the companion Delsarte-Genin
polynomials $\t P_n(x)$ as Christoffel transforms of the polynomials
$S_n(x)$. We show that the polynomials $Q_n(x)$ coincide with $\t
P_n(x)$.

In Section 6, we consider the Jacobi OPUC and show that the
corresponding polynomials $Q_n(x)$ coincide with the little -1
Jacobi polynomials.

In Section 7, a one parameter generalization of the previous
scheme - namely the Schur-Delsarte-Genin (SDG) map -  is proposed.
Starting from the same unitary matrices $L$ and $M$ we construct
a 3-diagonal matrix $K(\lambda)=L+\lambda M$ that depends on a
parameter $\lambda$. The matrix $J$ corresponds to $\lambda=1$.
The matrix $K(\lambda)$ generates a new family of orthogonal
polynomials $Q_n(x;\lambda)$ depending on the parameter $\lambda$.
A linear operator pencil is associated with the eigenvalue problem
for the matrix $K(\lambda)$.

In Section 8, we consider the special case
$a_0=a_1=\dots=0$ that corresponds to the OPUC $\Phi_n(z)=z^n$. The
polynomials $Q_n(x;\lambda)$ are then shown to be orthogonal on the
union of two intervals of the real axis.

In order to recover the orthogonality measure for the general
case, we need a special map from symmetric to non-symmetric
orthogonal polynomials that has been proposed by Chihara. This leads to a
nontrivial linear operator pencil $J_1 + \lambda J_2$. This
technique is described in Sections 8 and 9.

In Section 10 we relate the polynomials $Q_n(x;\lambda)$ to
special polynomials $P_n(x)$ independent of $\lambda$ and having
a known orthogonality measure.

Finally in Section 11, we apply this technique to Jacobi OPUC.
The corresponding polynomials $Q_n(x;\lambda)$ are shown to coincide
with the big -1 Jacobi polynomials which are orthogonal on the union of two
intervals on the real axis.

\section{Polynomials orthogonal on the unit circle and CMV matrices}
\setcounter{equation}{0}
Monic polynomials $\Phi_n(z)=z^n + O(z^{n-1})$ on the unit circle are defined through the recurrence relation
\be
\Phi_{n+1}(z) = z\Phi_n(z) - \bar a_{n} \Phi_n^*(z),
\lab{Sz_rec} \ee
with initial condition $\Phi_0(z) =1$
where
$$
\Phi_n^*(z) = z^n \overline{\Phi_n(1/\bar z)}.
$$
The recurrence coefficients $a_n=-\overline{\Phi_{n+1}(0)}$ are called
reflection parameters (sometimes also referred as the Schur, Geronimus,
Verblunsky, ... parameters). It is well known that the
reflection parameters are characterized by the property
$|a_n|<1$ for $n=0,1,2,\dots$.

{\bf In what follows we will consider only the case where $a_n$ are real}.

For convenience, we set
\be a_{-1}=-1 \lab{a_-1}.
\ee

The dual recurrence relation is
\be
\Phi_{n+1}^*(z) = \Phi_n^*(z) - a_n z \Phi_n(z) \lab{dual_rec_Phi}. \ee

It is standard to introduce the parameters
$$
r_n = \sqrt{1-a_n^2}
$$
(the arithmetic meaning of the square root is assumed).

Following \cite{Watkins} and \cite{CMV}, let us introduce the block-diagonal matrix
\vspace{5mm}
\be
 L =
 \begin{pmatrix}
  a_{0} & r_{0} &  &    \\
  r_{0} & -a_{0} &  &   \\
   &  &              a_{2} & r_{2}  \\
   &  &              r_{2} & -a_{2}  \\
  &  &    & &          a_{4} & r_{4}  \\
   &  &    & &         r_{4} & - a_{4}  \\
&   &  & & & & \ddots  \\
 \end{pmatrix}
\lab{L_def} \ee
and the block-diagonal matrix
\vspace{5mm}
\be
 M =
 \begin{pmatrix}
1 \\
& a_{1} & r_{1} &  &    \\
  &r_{1} & -a_{1} &  &   \\
   & &  &               a_{3} & r_{3}  \\
   & & &              r_{3} & -a_{3}  \\
  &  & &   & &           a_{5} & r_{5}  \\
   &  & &   & &         r_{5} & - a_{5}  \\
&   &  & & & & & \ddots  \\
 \end{pmatrix}.
\lab{M_def} \ee
Both $L$ and $M$ are block-diagonal unitary matrices. Define also the 5-diagonal unitary matrix $U$ \cite{Watkins}, \cite{CMV}
\be
U=LM \lab{def_U} \ee
(which is called the "Zigzag matrix" in \cite{Watkins} and the "CMV matrix" in \cite{Simon}).

In \cite{KN} it was proposed to consider the Hermitian matrix
\be
H = U + U^{\dagger}. \lab{HUU}
\ee
The authors of \cite{KN}
connected this matrix to a direct sum of two Jacobi matrices and
then arrived at the famous Szeg\H{o} relations between polynomials
orthogonal on the unit circle and on the real line. In the next section
we propose a slightly different approach to the matrix $H$
which connects the polynomials on the unit circle with a more
simple class of polynomials orthogonal on the interval. These
polynomials are related to symmetric orthogonal polynomials on
the same interval by a simple Christoffel transform.

\section{From CMV matrices to polynomials orthogonal on the real line}
\setcounter{equation}{0} On the natural basis $e_n, \;
n=0,1,2,\dots$ the matrix $H$ is 5-diagonal and Hermitian \ba
&&H e_n = r_n r_{n+1} e_{n+2} + r_n( a_{n+1} -a_{n-1}) e_{n+1} -2  a_{n} a_{n-1} e_n + \nonumber \\
&&r_{n-1}(a_{n} -a_{n-2}) e_{n-1} + r_{n-1} r_{n-2} e_{n-2}. \lab{H_e} \ea
Introduce the 3-diagonal matrix $J$ which is defined by
\be
J e_n = r_n e_{n+1} + (a_n - a_{n-1}) e_n + r_{n-1} e_{n-1}. \lab{J_e}
\ee
It is easily verified that
\be
H = J^2-2 \lab{H_J}.
\ee
Thus among the eigenvectors of the matrix $H$ there are the eigenvectors of the matrix $J$.
Moreover, it is seen that
\be J=L+M. \lab{J_LM}
\ee
Note that the matrices $L$ and $M$ are both unitary and satisfy the obvious
property $L^2=M^2=1$. Hence
\begin{equation}
J^2=(L+M)^2=L^2+M^2 +LM+ML = 2+U+U^*=2+H,
\end{equation}
thus recovering relation~\re{H_J}.

The matrix $J$ is a symmetric Jacobi matrix. It generates
orthonormal polynomials $\hat Q_n(x)$ satisfying the recurrence
relation
\begin{equation}
r_{n+1} \hat Q_{n+1}(x) + (a_n-a_{n-1}) \hat Q_n(x) + r_n \hat
Q_{n-1}(x) = x \hat Q_n(x).
\end{equation}
The corresponding monic orthogonal polynomials
\begin{equation}
Q_n(x)=\frac{\hat Q_n(x)}{r_1 r_2 \dots r_n} =x^n + O(x^{n-1})
\end{equation}
satisfy the recurrence
relation \be Q_{n+1}(x) + b_n Q_n(x) + u_n Q_{n-1}(x) = x Q_n(x)
\lab{rec_Q} \ee with initial conditions
$$
Q_0=1, \quad Q_1(x) =x-b_0,
$$
where the recurrence coefficients are
\be u_n = r_{n-1}^2 = 1-a_{n-1}^2, \quad b_n = a_n - a_{n-1}. \lab{rec_ub_Q} \ee
Note that
$$
u_0 =0, \quad b_0 = a_0+1
$$
due to the assumption $a_{-1}=-1$.

The relation between $z$ and $x$ is given by
\be
x=z^{1/2} + z^{-1/2}. \lab{x_z} \ee
Here it is assumed that one branch of the two-valued analytic function $z^{1/2}$ is chosen. In what follows we shall take
 $z^{1/2}= r^{1/2} e^{i\phi/2}$ when $z= r e^{i \phi}, \; 0 \le \phi \le 2\pi$.
Suppose that a point $z$ is moving counterclockwise on the unit circle $|z|=1$ starting from the point $z=1$ and
returning to the same point after a full rotation. Then the corresponding image point $x=z^{1/2} + z^{-1/2}$
is moving in one direction on the interval $[-2,2]$ starting from the point $2$ and ending at the point $-2$.
Thus the mapping \re{x_z} is a one-to-one correspondence between the punctured unit circle (i.e. $z \ne 1$)
and the interval $[-2,2]$. The only exception is the point $z=1$: it corresponds to the two
endpoints $\pm 2$ of the interval.

Generic polynomials $\Phi_n(z)$ under condition $|a_n|<1$ (with complex coefficients $a_n$) are orthogonal on the unit circle
\be
\int_{0}^{2 \pi} \Phi_n(e^{i \theta}) \bar \Phi_m(e^{-i \theta}) d \nu(\theta) = 0, \quad n \ne m, \lab{ort_nu} \ee
where $\nu(\theta)$ is a nondecreasing function in the interval $[0, 2 \pi]$.

For the case of real parameters $a_n$ we will assume that $d \nu(\theta) = \rho(\theta) d \theta$, where $\rho(\theta)$ is a positive weight function:
\be
\int_{0}^{2 \pi} \Phi_n(e^{i \theta}) \Phi_m(e^{-i \theta}) \rho(\theta) d \theta = 0, \quad n \ne m. \lab{ort_Phi} \ee
The function $\rho(\theta)$ satisfy the symmetric property
\be
\rho(2 \pi- \theta) = \rho(\theta)
\lab{sym_rho} \ee
which is valid for real coefficients $a_n$ verifying  the condition $-1 < a_n <1$ \cite{Simon}.

To end this section, note that $Q_n(x)$ are polynomials orthogonal with respect to a measure.
In fact, this measure can be found by means of the map proposed by Delsarte and Genin.

\section{The Delsarte-Genin map to symmetric polynomials on the interval}
\setcounter{equation}{0} Consider the symmetric polynomials
$S_n(x)$ satisfying the recurrence relation \be S_{n+1}(x) + v_n
S_{n-1}(x) = x S_n(x), \lab{rec_S} \ee where \be v_n =
(1+a_{n-1})(1-a_{n-2}) \lab{v_a} \ee Delsarte and Genin showed
\cite{DG1} that the polynomials $S_n(x)$ are related to the
polynomials $\Phi_n(z)$ as follows (see also \cite{Zhe1}) \be
S_n(x) = \frac{z^{-n/2}(\Phi_n(z) + \Phi_n^*(z))}{1-a_{n-1}}.
\lab{DG_map_SP} \ee The reciprocal map is \ba
&&\Phi_n(z) = \frac{z^{n/2}(z^{1/2} S_{n+1}(x) - \sigma_n S_n(x))}{z-1}, \nonumber \\
&&\Phi_n^*(z) = \frac{z^{(n+1)/2}(z^{1/2} S_{n+1}(x) - \sigma_n z^{1/2} S_n(x))}{1-z}   \lab{DG_map_PS} \ea
where
$$
\sigma_n = \frac{S_{n+1}(2)}{S_n(2)}.
$$
The reflection parameters $a_n$ are expressed as \be a_{n-1} =
1-\sigma_n \lab{a_A}. \ee The symmetric orthogonal polynomials
$S_n(x)$ are orthogonal on the interval $[-2,2]$: \be \int_{-2}^2
S_n(x) S_m(x) w(x) dx = 0, \quad n \ne m, \lab{ort_S_n} \ee where
the weight function $w(x)$ is symmetric $w(-x) = w(x)$ and is
related to the weight function $\rho(\theta)$ in the orthogonality
relation \re{ort_Phi} as follows~\cite{DG1,Zhe1} \be w(x) =
\frac{\rho(\theta)}{\sin(\theta/2)}, \quad x=2 \cos(\theta/2).
\lab{w_S} \ee Note that $\sin(\theta/2)>0$ for $0<\theta <2\pi$
and that the symmetry property \re{sym_rho} is equivalent to the
symmetry property $w(-x)=w(x)$.

There is an obvious "scaling" generalization of the Delsarte-Genin
mapping \cite{Zhe1}. For an arbitrary positive parameter $g$
define the map \be S_n(x;g) = \frac{g^n z^{-n/2} (\Phi_n(z) +
\Phi_n^*(z) )}{1-a_{n-1}}, \lab{SPP_g} \ee where the relation between
$x$ and $z$ is \be x= g(z^{1/2}+z^{-1/2}). \lab{xz_g} \ee It is
easily seen that the polynomials $S_n(x;g)$ are symmetric orthogonal
polynomials satisfying the recurrence relation
\be S_n(x;g) + g^2 v_n
S_{n-1}(x;g) = x S_n(x) \lab{rec_Sg} \ee
 with $v_n$ also given by~\re{v_a}. The polynomials $S_n(x)$ introduced in
\re{DG_map_SP} correspond to the case $g=1$. From \re{rec_Sg} it
is clear that \be S_n(x;g) = g^n S_n(x/g) \lab{SS_g} \ee i.e.
the polynomials $S_n(x;g)$ are obtained from the polynomials $S_n(x)$ by a
scaling transformation of the argument $x$. In particular, if
the polynomials $S_n(x)$ were orthogonal on the interval $[-2,2]$,
the polynomials $S_n(x;g)$ are orthogonal on the interval
$[2g,-2g]$.

The reciprocal transformation is given by the formulas \be \Phi_n(z) =
g^{-n-1} \frac{z^{n/2}(z^{1/2} S_{n+1}(x;g) - \sigma_n(g)
S_n(x;g))}{z-1}, \lab{Phi_Sg} \ee where
$$
\sigma_n(g) = \frac{S_{n+1}(2g;g)}{S_n(2g;g)} = g \sigma_n,
$$
and $$\sigma_n = \sigma_n(1)= \frac{S_{n+1}(2)}{S_n(2)}.$$

\section{The companion Delsarte-Genin map and orthogonality measures for SDG-transformed
polynomials}

In this section we find the orthogonality measure for $Q_n(x)$
by performing a Christoffel transformation of the polynomials from the previous section.
The resulting map will be called the Schur-Delsarte-Genin map (or, shortly, the SDG-map).

First, consider the polynomials $\tilde P_n(x)$ which are obtained
from $S_n(x)$ by the Christoffel transform \be \t P_n(x) =
\frac{S_{n+1} - A_n S_n(x)}{x+2},  \lab{tP_CT} \ee where \be A_n =
\frac{S_{n+1}(-2)}{S_n(-2)} = -\sigma_n =a_{n-1} -1
\lab{A_n_sigma} \ee (the last formula follows from the property
$S_n(-x)=(-1)^n S_n(x)$ for the symmetric orthogonal polynomials).

The polynomials $\t P_n(x)$ satisfy the recurrence relation
\be
\t P_{n+1}(x) + \t b_n \t P_n(x) + \t u_n \t P_{n-1}(x) = \t P_n(x), \lab{rec_t_P} \ee
where (see, e.g. \cite{ZhS}) the new recurrence coefficients are
\be
\t u_n = v_n \frac{A_n}{A_{n-1}}, \quad \t b_n = A_{n+1}-A_n. \lab{rec_t_ub} \ee
From \re{rec_t_ub} and \re{a_A}, we immediately see that the recurrence coefficients $\t b_n, \t u_n$ coincide with \re{rec_ub_Q}:
$$
\t u_n = 1-a_{n-1}^2, \quad \t b_n = a_n - a_{n-1}.
$$
This means that the polynomials $\t P_n(x)$ and $Q_n(x)$ coincide. Thus the (non-symmetric) polynomials
$Q_n(x)$ are obtained from the symmetric polynomials $S_n(x)$ by the Christoffel transform \re{tP_CT}.
This allows expressing the polynomials $Q_n(x)$ in terms of the circle polynomials $\Phi_n(z)$:
\be
Q_n(x) = z^{-n/2} \: \frac{\Phi_n^*(z) + z^{1/2} \Phi_n(z)}{1+z^{1/2}}. \lab{Q_Phi} \ee
The reciprocal transformation from the polynomials $S_n(x)$ to the polynomials $Q_n(x)$
is given by the Geronimus transform
\be
S_n(x)=Q_n(x) -B_n Q_{n-1}(x),
\lab{S_Q} \ee
where
\be
B_n=2-A_n =1+a_{n-1}.
\ee
The polynomials $Q_n(x)$ are orthogonal on the interval $[-2,2]$:
\be
\int_{-2}^2 Q_n(x) Q_m(x) w(x) (x+2) dx =0, \quad n \ne m, \lab{ort_Q} \ee
where the function $w(x)$ is given by \re{w_S}. Formula \re{ort_Q} follows from the definition of the Christoffel transform \cite{Sz}.

Consider the special value $Q_n(2)$. From \re{Q_Phi} we have ($x=2$ corresponds to $z=1$)
\be
Q_n(2)=\Phi_n(1) = (1-a_0)(1-a_1)\dots(1-a_{n-1}).
\ee
Using these formulas, one can easily obtain an expression for the polynomials $\Phi_n(z)$ in terms of the
polynomials $Q_n(x)$:
\be
\Phi_n(z) = z^{n/2} \:\frac{(z^{-1/2}-a_n)Q_n(x) -Q_{n+1}(x)}{1-z^{1/2}}. \lab{Phi_Q} \ee
Note that apart from the polynomials $Q_n(x)$ one can introduce the adjacent polynomials $Q_n^{(-)}(x)$
which are obtained from $S_n(x)$ by the Christoffel transformation which is adjacent with respect to \re{tP_CT}:
\be
Q_n^{(-)}(x) = \frac{S_{n+1} - A_n^{(-)} S_n(x)}{x-2}, \quad A_n^{(-)} =
\frac{S_{n+1}(2)}{S_n(2)} = - A_n = 1-a_{n-1}. \lab{Q-_CT} \ee
The monic polynomials $Q_n^{(-)}(x)$ satisfy the recurrence relation
\be
Q_{n+1}^{(-)}(x) + b_n^{(-)} Q_n^{(-)}(x) + u_n^{(-1)} Q_{n-1}^{(-)}(x) = x Q_n^{(-)}(x), \lab{rec_Q-} \ee
where
\be
b_n^{(-)}=a_{n-1}-a_n = -b_n, \quad u_n^{(-1)} = 1-a_{n-1}^2 = u_n. \lab{bu_minus} \ee
Comparing \re{rec_Q-} with \re{rec_Q}, we see that the polynomials $Q_n^{(-)}(x)$ have the same
recurrence coefficients $u_n$ as the polynomials $Q_n(x)$ but with coefficient $b_n$ that have
the opposite sign. This means that these polynomials almost coincide:
\be
Q_n^{(-)}(x) = (-1)^n Q_n(-x) \lab{Q-Q}. \ee
The polynomials $Q_n^{(-)}(x)$ are orthogonal on the interval $[-2,2]$
\be
\int_{-2}^2 Q_n^{(-)}(x) Q_m^{(-)}(x) w(x) (x-2) dx =0, \quad n \ne m. \lab{ort_Q_minus} \ee
We see that the polynomials $Q_n(x)$ are more convenient to construct the polynomials $\Phi_n(z)$ than the
polynomials $S_n(x)$. Indeed, the reflection parameters $a_n$ are immediately determined, up to a sign, from the recurrence coefficient $u_n=1-a_{n-1}^2$:
\be
a_{n} =\pm \sqrt{1-u_{n+1}}. \lab{a_u} \ee
Choosing the sign in \re{a_u} selects the polynomials $Q_n(x)$ or the polynomials $Q_n^{(-)}(x)$.

We can consider a third set of orthogonal polynomials $T_n(x)$ which is obtained from
the polynomials $Q_n(x)$ by the Christoffel transform at $x=2$:
\be
T_n(x)= \frac{Q_{n+1}(x) - \frac{Q_{n+1}(2) }{Q_{n}(2) }Q_n(x)}{x-2} =  \frac{Q_{n+1}(x) - (1-a_{n})Q_n(x)}{x-2}.
\lab{T_Q} \ee
Simple calculations show that the polynomials $T_n(x)$ are symmetric orthogonal polynomials satisfying the recurrence relation
\be
T_{n+1}(x) + w_n T_{n-1}(x) = xT_n(x), \quad T_0=1, \; T_1 =x, \lab{T_rec} \ee
where
\be
w_n=(1+a_{n-1})(1-a_n). \lab{w_rec} \ee
Note that the polynomials $T_n(x)$ can be obtained from the polynomials $S_n(x)$ by a "symmetric Christoffel transform"
(see, e.g. \cite{Sz}, \cite{SZ})
\begin{equation}\lab{T_CT_S}
\begin{split}
T_n &= \frac{S_{n+2}(x) - \frac{S_{n+2}(2)}{S_n(2)}  S_n(x)}{x^2-4}= \frac{S_{n+2}(x) - (1-a_n)(1-a_{n-1})  S_n(x)}{x^2-4}\\
&=\frac{xS_{n+1}(x) -2(1-a_{n-1})S_n(x)}{x^2-4}.
\end{split}
\end{equation}
Using formula \re{DG_map_SP} and the recurrence relations for the polynomials $\Phi_n(z)$
we can express the polynomials $T_n(x)$ in terms of the polynomials $\Phi_n(z)$:
\be
T_n(x) = \frac{z^{-n/2} (z\Phi_n(z) - \Phi^*_n(z)) }{z-1}.  \lab{T_Phi} \ee
The polynomials $T_n(x)$ are orthogonal on the interval $[-2,2]$:
\be
\int_{-2}^2 T_n(x) T_m(x) (4-x^2) w(x) dx =0, \quad n \ne m. \lab{ort_T} \ee
Thus with the polynomials $\Phi_n(z)$ orthogonal on the unit circle with real reflection parameters $|a_n|<2$
one can relate 3 types of polynomials $S_n(x), Q_n(x), T_n(x)$, orthogonal on the interval $[-2,2]$ of the real line.
Two of these polynomials, $S_n(x)$ and $T_n(x)$, are symmetric, whereas the polynomials $Q_n(x)$ are non-symmetric.

Note that the polynomials $S_n(x)$ and $T_n(x)$ coincide with
the so-called "companion polynomials" considered in \cite{BR}.
Formulas \re{DG_map_PS} express the polynomials $\Phi_n(z)$ in terms
of the polynomials $S_n(x)$. Note also that the relations between
the symmetric polynomials $T_n(x)$ and the non-symmetric polynomials
$Q_n(x)$ (without the identification of the connection to OPUCs) were studied by L.M.~Chihara and
T.S.~Chihara in~\cite{ChCh}.

\section{From the Jacobi OPUC to the little -1 Jacobi polynomials}
\setcounter{equation}{0}

The unit circle analogues of the Jacobi polynomials were proposed
by Szeg\H{o}~\cite{Sz} and studied by Golinskii~\cite{Gol} and
Badkov~\cite{Badkov83},\cite{Badkov}. These Jacobi OPUC
$\Phi_n(e^{i\theta})$ are orthogonal on the unit circle with
respect to the weight function
\begin{equation}
\rho(\theta)=(1-\cos\theta)^{\xi+1/2}(1+\cos\theta)^{\eta+1/2},
\end{equation}
where $\theta \in[0,2\pi]$. It is assumed that $\xi>-1, \;
\eta>-1$.

Their parameters $a_n$ have the following expression \cite{Badkov83} \be a_n
= \left\{ \frac{\eta-\xi}{n+\xi+\eta+2} \quad \mbox{if} \quad n
\quad \mbox{is} \quad \mbox{even} \atop
-\frac{1+\xi+\eta}{n+\xi+\eta+2} \quad \mbox{if} \quad n \quad
\mbox{is} \quad \mbox{odd}  \right . .\lab{gen_ultra_a} \ee The
condition $a_{-1}=-1$ obviously holds.

Note that the Jacobi OPUC are connected with the generalized
Gegenbauer (ultraspherical) polynomials through the DG map.

Indeed, the generalized Gegenbauer polynomials \cite{Chi},
\cite{Bel} are symmetric orthogonal polynomials satisfying
\re{rec_S} with \be v_n = \left\{ \frac{n(n+2
\xi)}{(n+\xi+\eta)(n+\xi+\eta+1)} \quad \mbox{if} \quad n \quad
\mbox{is} \quad \mbox{even} \atop \frac{(n+2\eta+1)(n+2 \xi+2
\eta+1)}{(n+\xi+\eta)(n+\xi+\eta+1)} \quad \mbox{if} \quad n \quad
\mbox{is} \quad \mbox{odd}  \right . ,\lab{gen_ultra_v} \ee where
$\xi>-1, \: \eta>-1$. These polynomials are orthogonal on the
interval $[-2,2]$: \be \int_{-2}^2 S_n(x) S_m(x) (4-x^2)^{\xi}
|x|^{2\eta+1} dx = 0, \quad n \ne m. \lab{ort_gen_ultra} \ee Using
the Delsarte-Genin map \re{DG_map_PS} and formula \re{v_a}, we can
easily verify that the Jacobi OPUC  $\Phi_n(e^{i\theta})$
correspond to the generalized Gegenbauer polynomials $S_n(x)$.

Consider now the companion polynomials $Q_n(x)$ corresponding to
the parameters \re{gen_ultra_a}. They are orthogonal on the
interval $[-2,2]$ with respect to the function
$$
w(x) = (x+2) (4-x^2)^{\xi} |x|^{2\eta+1}.
$$
The recurrence coefficients $b_n$, $u_n$ for these polynomials are
explicitly obtained from formulas \re{rec_ub_Q}, where the
coefficients $a_n$ are given by \re{gen_ultra_a}.

Choose the parametrization
$$
\xi=\frac{\alpha-1}{2}, \quad \eta=\frac{\beta-1}{2}.
$$
A simple analysis shows that the polynomials $Q_n(x)$ coincide (up
to an obvious affine transformation  of the argument) with the little
-1 Jacobi polynomials $P_n^{(-1)}(x;\alpha,\beta)$ introduced in
\cite{VZ_little}. The little -1 Jacobi polynomials possess a
remarkable property: they are polynomial  eigenfunctions of a
Dunkl-type differential operator of the first order. Hence these
polynomials provide a new nontrivial example of "classical"
orthogonal polynomials on the real line (for details see
\cite{VZ_little}.

Note that similar polynomials were considered in \cite{ChCh},
where the recurrence coefficients $b_n, u_n$ were calculated
explicitly. These polynomials belong to a special class of
orthogonal polynomials and were studied by another approach
in \cite{AAR}. Our approach yields explicit expressions for
the recurrence coefficients and the orthogonal polynomials $Q_n(x)$ by a
straightforward  application of the correspondence formulas between polynomials
orthogonal on the unit circle and on the interval.

\section{A one-parameter generalization of the companion polynomials}
\label{reductionQ}
Consider a linear pencil of matrices
\be
K(\lambda) = L + \lambda M \lab{lin_K} \ee
with an arbitrary parameter $\lambda$.

We have \be K^2(\lambda) = L^2 + \lambda^2 M^2 + \lambda (LM + ML)
= 1+\lambda^2 + \lambda (J^2-2) = 1+ \lambda^2 + \lambda H
\lab{K_J_H}. \ee We see that the Hermitian matrix $H$ can be
expressed in terms of $K^2(\lambda)$ for every $\lambda$.  On the
other hand, the matrix $K(\lambda)$ is a nondegenerate symmetric
tri-diagonal matrix
\[
 K(\lambda) =
 \begin{pmatrix}
a_0 + \lambda & r_0 \\
r_0 & -a_{0}+\lambda a_1 & \lambda r_{1} &  &    \\
  &\lambda r_{1} & a_2-\lambda a_{1} &  r_2 &   \\
   & &  r_2 & -a_2+\lambda a_3 &   \lambda r_3 &  \\
  &   &  & &  \ddots  \\
 \end{pmatrix}.
\]
Here, nondegenerate means that all off-diagonal entries of the matrix $K(\lambda)$ are nonzero
if $\lambda \ne 0$ and $|a_i| \ne 1$.

Hence one can define a family of formal monic orthogonal
polynomials $Q_n(x;\lambda)$ depending on the argument $x$ and an
additional parameter $\lambda$. These polynomials are uniquely defined
 through the 3-term recurrence relation
\be
Q_{n+1}(x;\lambda) + b_n(\lambda) Q_n(x;\lambda) + u_n(\lambda)
Q_{n-1}(x;\lambda) = x Q_n(x;\lambda), \lab{rec_Q_lambda}
\ee
where
\be b_n(\lambda) = \left\{a_n - \lambda a_{n-1} \quad \mbox{if}
\quad n \quad \mbox{is} \quad \mbox{even}  \atop \lambda a_n -
a_{n-1} \quad \mbox{if} \quad n \quad \mbox{is} \quad \mbox{odd}
\right . ,\lab{b_lambda} \ee and \be u_n(\lambda) =
\left\{\lambda^2 (1-a_{n-1}^2) \quad \mbox{if} \quad n \quad
\mbox{is} \quad \mbox{even} \atop  1-a_{n-1}^2 \quad \mbox{if}
\quad n \quad \mbox{is} \quad \mbox{odd}  \right . .\lab{u_lambda}
\ee Note that the eigenvalue problem for the orthogonal polynomials
$Q_n(x;\lambda)$ can be presented in algebraic form as \be (L +
\lambda M -x I) \vec q =0, \lab{lambda_x_q} \ee where $\vec q$ is a
vector constructed from the (non-monic) polynomials $Q_n(x;\lambda)$.
Equation \re{lambda_x_q} contains two parameters $\lambda$ and $x$
and belongs to the class of the so-called multi-parameter eigenvalue problems
\cite{Atkinson}, \cite{Sleeman}.

We already considered the case $\lambda=1$ that leads to the orthogonal polynomials $Q_n(x)$ with a positive measure
on the interval $[-2,2]$. For general complex values of the parameter $\lambda$ it is still
unclear what the spectral properties of the corresponding orthogonal polynomials are
because the matrix $K(\lambda)$ is non-Hermitian.
Nevertheless, when $\lambda$ is real, the polynomials
$Q_n(x;\lambda)$ are orthogonal on the real axis with a positive measure.
Indeed, as seen from \re{u_lambda} the recurrence coefficients $u_n$ are positive for all $n>0$
and for all real values of $\lambda$.  In the special case $\lambda=-1$, there is a simple relation
between the polynomials $Q_n(x;-1)$ and the polynomials $(-1)^n\Phi_n(-z)$ on the unit circle.

Indeed, in the case $\lambda=-1$, the recurrence coefficients of the polynomials $Q_n(x;-1)$ become
\be
b_n(-1)= (-1)^n (a_n + a_{n-1}), \quad u_n(-1) = 1-a_{n-1}^2.
\lab{bu_-1} \ee
The coefficients \re{bu_-1} are obtained from the coefficients \re{bu_minus}
by a simple transformation $a_n \to (-1)^{n+1} a_n$ corresponding to
the monic polynomials $\tilde \Phi_n(z) = (-1)^n \Phi_n(-z)$ which
are also polynomials orthogonal on the unit circle (since obtained from a
simple reflection of the argument $z$ with respect to the real axis). The moments $\tilde c_n$ corresponding to the polynomials $\tilde \Phi_n(z)$ are $\tilde c_n =(-1)^n c_n$.

It is easily verified that
$$
A_n \equiv \frac{Q_{n+1}(-2;-1)}{Q_n(-2;-1)} = -(1+(-1)^n a_n)
$$
and hence, one can introduce new symmetric polynomials $\tilde T_n(x)$ obtained from $Q_n(x;-1)$ by the Christoffel transform
$$
\tilde T_n(x) = \frac{Q_{n+1}(x;-1) - A_n Q_n(x;-1)}{x+2}.
$$
These polynomials satisfy the recurrence relation
$$
\tilde T_{n+1}(x) + \tilde v_n T_{n-1}(x) = x \tilde T_n(x)
$$
with
$$
\tilde v_n = (1+(-1)^n a_n)(1+(-1)^n a_{n-1}).
$$

Note that Watkins proposed \cite{Watkins} to consider the unitary linear pencil problem
\be
(L + \lambda M) \vec v =0 \lab{eigen_v} \ee
which leads to Szeg\H{o} polynomials orthogonal on the unit circle.
In this case $x=0$ and the parameter $\lambda$ belongs to the unit circle $|\lambda|=1$.

\section{The 1-periodic case}

In order to better understand  the nature of the polynomials
$Q_n(x;\lambda)$ let us consider the simplest periodic case, that
is, $a_0=a_1=\dots=0$ implying that $r_0=r_1=\dots=1$. The
corresponding OPUC are simple monomials \cite{Simon} \be \Phi_n(z)
= z^n . \lab{monom_Phi} \ee

The symmetric polynomials $S_n(x)$ obtained by the DG transform
\re{DG_map_SP} coincide with the Chebyshev polynomials of the
first kind \be S_0=1, \quad S_n(x) = 2 \cos(n \theta/2), \quad
n=1,2,\dots, \quad x= 2 \cos(\theta/2). \lab{S_Cheb} \ee They are
orthogonal on the interval $[-2,2]$ with the weight function \be
w_0(x) = \frac{1}{\sqrt{4-x^2}} . \lab{CH_w} \ee In order to recover
the orthogonality measure for the polynomials $Q_n(x, \lambda)$ we
need some additional tools. Namely, let us introduce the following function
\begin{equation}
    m(z,\lambda)=\langle (z-K(\lambda))^{-1}e_0,e_0\rangle ,\quad z\in\dC_+,
\end{equation}
where $e_0=(1,0,\dots)^\top$. It is called the $m$-function (or Weyl function) of
$K(\lambda)$~\cite[Section 1.2]{Simon}. It is a standard fact that
\begin{equation}\label{Ric}
m(z,\lambda)=-\frac{1}{z-\lambda+m_1(z,\lambda)},\quad \lambda\in{\mathbb R}_+,
\end{equation}
where $m_1(z,\lambda)$ is the $m$-function of the matrix
\be
 K_{1}(\lambda)=
 \begin{pmatrix}
 0 & \lambda &  &    \\
 \lambda  & 0 &  1 &   \\
  &  1 & 0 &   \lambda  &  \\
     &  & \lambda &  \ddots  \\
     &&&\\
 \end{pmatrix},\quad \lambda\in{\mathbb R}_+.
\ee
Let us also consider the following family of matrices
\be
K_{per}(\lambda)=
 \begin{pmatrix}
 0 & 1 &  &    \\
 1 & 0 &  \lambda &   \\
  &  \lambda & 0 &   1  &  \\
     &  & 1 &  \ddots  \\
     &&&\\
\end{pmatrix}=
K(\lambda)-\text{diag}(\lambda,0,0,\dots),\quad \lambda\in{\mathbb R}_+.
\ee
Due to~\eqref{Ric}, we have that
\begin{equation}\label{relat_per}
m(z,\lambda)=\frac{m_{per}(z,\lambda)}{1+\lambda m_{per}(z,\lambda)}.
\end{equation}

 Next, observe that $m_{per}$ satisfies the following equation
\be
m_{per}(z)=-\frac{1}
{z-{ \displaystyle \frac{1} {z+\lambda^2m_{per}(z,\lambda)} } },
\ee
which is equivalent to the algebraic equation
\be
\lambda^2zm_{per}^2+
\left(z^2-1+\lambda^2\right)m_{per}+z=0.
\ee
Therefore, we have that
\begin{equation}\label{weyl_per}
m_{per}(z,\lambda)=\frac{1}{2z\lambda^2}\left(-z^2-\lambda^2+1+
\sqrt{(z^2-\lambda^2-1)^2-4\lambda^2}\right),
\end{equation}
where the branch of the square root is chosen
such that $m_{per}(z,\lambda)\sim-\frac{1}{z}$ as $z\to\infty$.
Recall that the spectrum $\sigma(K)$ of the operator $K$ is the set of all complex numbers $\zeta$ for which
$K-\zeta I$  is not invertible. The essential spectrum $\sigma_{ess}(K)$ is defined as
$\sigma_{ess}(K)=\sigma(K)\setminus\sigma_{d}(K)$, where $\sigma_{d}(K)$ is the set of all isolated
eigenvalues with finite multiplicity.

The Floquet theory for semi-finite periodic self-adjoint Jacobi matrices~\cite{Ger57} says
that the essential spectrum of $K_{per}(\lambda)$ admits the representation
\be\label{Es_sp}\begin{split}
\sigma_{ess}(K_{per}(\lambda))=\{z\in\dR: (z^2-(1+\lambda^2)/\lambda\in[-2,2])\}=\\
=[-\lambda-1,-|\lambda-1|]\cup[|\lambda-1|,\lambda+1].
\end{split}
\ee
Moreover, $K_{per}(\lambda)$ has 0 as an eigenvalue for $\lambda>1$ ~\cite{Ger57}.
So, according to Weyl's theorem we have
that $\sigma_{ess}(K(\lambda))=\sigma_{ess}(K_{per}(\lambda))$ since $K-K_{per}$ is a self-adjoint one-dimensional
operator. Furthermore,  $\sigma(K(\lambda))=\sigma_{ess}(K(\lambda))$ because $m(z,\lambda)$
does not have isolated poles.
Thus, from~\eqref{Es_sp} we clearly see what the $\lambda$-dynamics
 of  the spectrum of $K(\lambda)$ is.
Namely, the spectrum starts with two points 1 and -1 at $\lambda=0$. Then, the points become
 two different intervals moving towards each other while $\lambda$ runs from 0 to 1. At the time $\lambda=1$,
 the intervals touch each other. After that, they are moving towards infinity in the opposite directions as $\lambda\to\infty$.

Finally,  substituting~\eqref{weyl_per} into~\eqref{relat_per} and applying the Stieltjes--Perron inversion formula,
we arrive to the fact that the polynomials $Q_n(z,\lambda)$ are orthogonal with respect to the weight function~\cite{Ch57}
\be
w(t)=\begin{cases}
\frac{\sqrt{4\lambda^2-(t^2-\lambda^2-1)^2}}{-\lambda (t^2-2\lambda^2t+\lambda^2-1)},
\quad t\in[|\lambda-1|,\lambda+1],\\
\frac{\sqrt{4\lambda^2-(t^2-\lambda^2-1)^2}}{\lambda (t^2-2\lambda^2t+\lambda^2-1)}, \quad t\in[-\lambda-1,-|\lambda-1|].
\end{cases}
\ee

In particular, for $\lambda=1$ we get that \be
w(t)=\frac{\sqrt{4-t^2}}{2-t}=\sqrt{\frac{2+t}{2-t}}=w_0(t)(t+2).
\lab{w_3} \ee From \re{w_3} it is seen that the polynomials
$Q_n(x)= Q_n(x;1)$ coincide with the Chebyshev polynomials of the
third kind \cite{MH}. Similarly, for $\lambda=-1$ the polynomials
$Q_n(x;-1)$ coincide with the Chebyshev polynomials of the fourth
kind.

In fact, the behavior of the set of orthogonality is typical for the entire class of general matrices $K(\lambda)$.
Moreover, the weight functions have a similar structure as well.
However, before dealing with the general case we have to develop a special technique.

\section{From the Chihara map to the generic SDG-maps}

Here, we recall some well-known results in the theory of orthogonal polynomials.
The main goal of this section is to elaborate the machinery for
the SDG-maps in the case of arbitrary real $\lambda$, i.e.
the generic SDG-maps.

Let us start with the following simple lemma
\begin{lemma}\label{lemma1}
Let $P_n(x)$ be monic orthogonal polynomials satisfying the
recurrence relation
\be P_{n+1}(x) +(\theta - A_n - C_n)P_n(x) +
C_n A_{n-1} P_{n-1}(x) = xP_n(x), \lab{3-term_AC} \ee
with the standard initial conditions
\be
P_0=1, \; P_1(x) = x-\theta+A_0 .
\ee
Assume that the real coefficients $A_n,C_n$ are such that
$A_{n-1}>0,\; C_n
>0, \; n=1,2\dots$ and that $C_0=0$. Take $\theta$ to be an arbitrary real
parameter such that $P_n(\theta)\ne 0$ for $n=1,2,\dots$.

Define the new polynomials \be \t P_n(x) = \frac{P_{n+1}(x) - A_n
P_n(x)}{x-\theta}. \lab{tP_CT1} \ee Then the monic polynomials $\t
P_n(x)$ are orthogonal and satisfy the recurrence relation \be \t
P_{n+1}(x) +(\theta - A_n - C_{n+1}) \t P_n(x) + C_n A_{n} \t
P_{n-1}(x) = x \t P_n(x). \lab{3-term_tP} \ee The inverse
transformation from the polynomials $\t P_n(x)$ to the polynomials
$P_n(x)$ is given by the formula \be P_n(x) = \t P_n(x) - C_n \t
P_{n-1}(x) . \lab{GT} \ee
\end{lemma}

To prove Lemma~\ref{lemma1} it is sufficient to observe that relation~\eqref{3-term_AC}
provides us with the explicit $LU$-factorization of the monic Jacobi matrix $J-\theta$
corresponding to the polynomials $P_n(x)$. The rest is an obvious reformulation
of the theory of spectral (or Darboux) transformations \cite{ZhS},
\cite{BM}.

Note also that \be A_n = \frac{P_{n+1}(\theta)}{P_n(\theta)}
\lab{A_PP} \ee which can easily be verified directly from
\re{3-term_AC}.

\begin{remark}
If the polynomials $P_n(x)$ are orthogonal with respect to a
weight function $w(x)$ on the interval $[\alpha,\beta]$ then the
polynomials $\t P_n(x)$ are orthogonal with respect to the
weight function $(x-\theta)w(x)$ on the interval $[\alpha,\beta]$.
\end{remark}

The next Lemma will be useful in the identification  of polynomial
systems with known families of orthogonal polynomials.

\begin{lemma}[\cite{Chi_Bol}]\label{lemma2}
Let $P_n(x)$ and $\t P_n(x)$ be two systems of orthogonal
polynomials defined as in the previous Lemma. Also, assume that
the polynomials $P_n(x)$ are orthogonal with respect to a
weight function $w(x)$ on the finite interval $[\alpha,\beta]$ and take any
$\theta\in(-\infty,\alpha)$. Finally, define the following monic polynomials
\be S_{2n}(x) = P_n(x^2+\alpha-c^2) ,
\quad S_{2n+1}(x) = (x-\chi) \t P_n(x^2+\alpha-c^2), \lab{S_PP} \ee
where $c$ is a positive number such that $\sqrt{\alpha-\theta}\le c$
and $\chi$ is a real number defined by the relation $\theta=\chi^2+\alpha-c^2$.
Then the polynomials $S_n(x)$ are orthogonal with respect to the weight function
\be
\sign x(x+\chi)w(x^2+\alpha-c^2)
\ee
on the union of two intervals $[-\sqrt{\beta-\alpha+c^2},-c]\cup [c,\sqrt{\beta-\alpha+c^2}]$
and satisfy the recurrence
relation \be S_{n+1} + (-1)^n \chi S_n(x) + v_n S_{n-1}(x) =
xS_n(x), \lab{rec_S1} \ee
where $v_n$ are given by \be v_{2n} =-C_n, \; v_{2n+1}=-A_n.
\lab{v_CA} \ee
The converse statement is also true.
\end{lemma}

The orthogonality part of this Lemma can be checked by straightforward computations.
Then, substituting~\eqref{S_PP} into the corresponding three-term recurrence relations
we arrive at~\eqref{rec_S1} and~\eqref{v_CA}. For more details, see~\cite{Chi_Bol}.

In what follows, we shall call the polynomials $\t P_n(x)$, the companion
polynomials with respect to $P_n(x)$ (sometimes the polynomials
$\t P_n(x)$ are referred to as  the kernel polynomials
\cite{Chi}).

%Now, we are in a position to consider the polynomials $Q_n(x;\lambda)$ again. Namely,
%let us recall that we reduced them to the polynomials satisfying~\eqref{rec_tQ}
%with coefficients either~\eqref{tbu_lambda} or~\eqref{tbu_lambda_1}.

\section{A scheme for reducing the polynomials $Q_n(x;\lambda)$ to $\lambda$-independent orthogonal polynomials}

In this Section, we show how to reduce the polynomials $Q_n(x;\lambda)$ defined by~\eqref{rec_Q_lambda}
to a $\lambda$-independent system of orthogonal polynomials.

First, let us apply Lemma~\ref{lemma1} to $Q_n(x;\lambda)$ for $\theta=\lambda-1$ and  $\theta=\lambda+1$.
Namely, if we put $\theta=\lambda-1$ then we can deduce that
\be
A_n=\left\{-1-a_{n} \quad \mbox{if} \quad n \quad \mbox{is} \quad \mbox{even} \atop
-\lambda(1-a_{n}) \quad \mbox{if} \quad n \quad \mbox{is} \quad \mbox{odd}  \right .
\lab{A_n_lambda},\quad
C_n=\left\{ -\lambda(1+a_{n-1}) \quad \mbox{if} \quad n \quad \mbox{is} \quad \mbox{even} \atop
-1+a_{n-1}\quad \mbox{if} \quad n \quad \mbox{is} \quad \mbox{odd}  \right . .
 \ee
Next, introducing the polynomials
\be
\t Q_n(x;\lambda) = \frac{Q_{n+1}(x;\lambda) - A_n Q_n(x;\lambda)}{x-\lambda+1}
\lab{tQ_lambda} \ee
we arrive at the the recurrence relation
\be
\t Q_{n+1}(x;\lambda)  + \t b_n \t Q_{n}(x;\lambda) + \t u_n \t Q_{n-1}(x;\lambda)=x \t Q_{n}(x;\lambda), \lab{rec_tQ} \ee
where
\be
\t b_n = (-1)^n (\lambda+1), \quad \t u_n =-\lambda (1+a_n)(1+(-1)^n a_{n-1}).
\lab{tbu_lambda} \ee
%Note that $\t u_n>0$ if $\lambda<0$, hence for all values $\lambda<0$ the polynomials $\t Q_n(x;\lambda)$ have a positive measure on the real line.

Similarly, for $\theta=\lambda+1$ we have that
\be
A_n= \left\{1-a_{n} \quad \mbox{if} \quad n \quad \mbox{is} \quad \mbox{even} \atop
\lambda(1-a_{n}) \quad \mbox{if} \quad n \quad \mbox{is} \quad \mbox{odd}  \right .,\quad
C_n= \left\{ -\lambda(1+a_{n-1}) \quad \mbox{if} \quad n \quad \mbox{is} \quad \mbox{even} \atop
-1-a_{n-1}\quad \mbox{if} \quad n \quad \mbox{is} \quad \mbox{odd}  \right .
.
\lab{A_n_lambda_1} \ee
Then the corresponding transformed polynomials $\t Q_n(x;\lambda)$ have the recurrence coefficients
\be
\t b_n =(-1)^n(\lambda-1), \quad \t u_n = \lambda(1+a_{n-1})(1-a_n). \lab{tbu_lambda_1} \ee

Clearly,
both cases can be rewritten as follows
\begin{equation}\label{UniF}
\t Q_{n+1}(x;\lambda)  +(-1)^n(d_1\lambda+d_0) \t Q_{n}(x;\lambda) +
\lambda u_n^* \t Q_{n-1}(x;\lambda)=x \t Q_{n}(x;\lambda),
\end{equation}
where $d_0$, $d_1$, and $u_n^*>0$ do not depend on $\lambda$.

The following Lemma reduces the polynomials $\t Q_{n+1}(x;\lambda)$ to the polynomials obtained in
Lemma~\ref{lemma2}.

\begin{lemma}\label{lemma3} Let the polynomials $\t Q_{n}(x;\lambda)$ satisfy~\eqref{UniF} and be orthogonal
with respect to a weight function $w_\lambda(x)$ on the set $E_\lambda$ for positive $\lambda$.
Then the monic polynomials defined via
\be
S_n(x;\lambda)=(\sqrt{\lambda})^{-n}Q_n(\sqrt{\lambda}x;\lambda)
\ee
are orthogonal with respect to the weight function $w_\lambda(\sqrt{\lambda}x)$ on the set
$E_\lambda^*=\{x:\,\, \sqrt{\lambda}x\in E_\lambda\}$
and satisfy
\begin{equation}\label{UniF_1}
S_{n+1}(x;\lambda)  +(-1)^n\chi S_{n}(x;\lambda) +
u_n^* S_{n-1}(x;\lambda)=x S_{n}(x;\lambda),
\end{equation}
where $\chi=d_1\sqrt{\lambda}+d_0\frac{1}{\sqrt{\lambda}}$.

The converse is also true.
\end{lemma}

The proof is straightforward by making the substitution $x\to\sqrt{\lambda}x$.

Now, we can apply Lemma~\ref{lemma2} to represent the polynomials $S_n(x; \lambda)$ in terms of
a $\lambda$-independent system of orthogonal polynomials. Also, Lemma~\ref{lemma2} gives us the explicit
formula for the weight function of the polynomials $S_n(x; \lambda)$.

Summing up, let us notice that the scheme given in this Section gives us the possibility
to reduce the $\lambda$-dependent system of polynomials $Q_n(x;\lambda)$ to a system of $\lambda$-independent
orthogonal polynomials $P_n(x)$ given by~\eqref{3-term_AC}. Namely, we reduced $Q_n(x;\lambda)$ to the
polynomials $\t Q_n(x;\lambda)$ satisfying~\eqref{UniF}.
Then Lemma~\ref{lemma3} leads us to the polynomials $S_n$, which according to Lemma~\ref{lemma2} can be
constructed from the polynomials $P_n(x)$ given by~\eqref{3-term_AC}.

Finally, it is worth noting that Lemma~\ref{lemma2} and Lemma~\ref{lemma3} give an efficient way to
construct explicit families of Jacobi matrices of the form
\[
J_1+\lambda J_2,
\]
where $J_1$ and $J_2$ are bidiagonal matrices. The corresponding
eigenvalue problems with spectral parameter $x$ lead to
orthogonal polynomials and to Laurent biorthogonal polynomials if
we fix $x$ and consider the problem with respect to $\lambda$
\cite{HR}, \cite{ZheL}.

\section{From the Jacobi OPUC to the big -1 Jacobi polynomials}
\setcounter{equation}{0}

To fit the big -1 Jacobi polynomials in this context, let us consider
the polynomials
\be
Q_n^*(x; \lambda)=g^{-n}Q_n(gx;\lambda).
\ee
Then~\eqref{rec_Q_lambda} yields the following
\be
Q_{n+1}^*(x;\lambda) + b_n^*(\lambda) Q_n^*(x;\lambda) + u_n^*(\lambda)
Q_{n-1}^*(x;\lambda) = x Q_n^*(x;\lambda),
\ee
where the coefficients are as follows
\be
b_n^*=\frac{1}{g}(\lambda-1-A_n-C_n),\quad u_n^*=\frac{1}{g^2}C_nA_{n-1},
\ee
and $A_n$, $C_n$ are defined in~\eqref{A_n_lambda}.

Next we compute the coefficients $A_n$ and $C_n$ for the
reflection parameters $a_n$ given by \re{gen_ultra_a}. We have
\be
A_n=\left\{-\frac{n+2+2\eta}{n+2+\xi+\eta} \quad \mbox{if} \quad n \quad \mbox{is} \quad \mbox{even} \atop
-\lambda\frac{n+3+2\xi+2\eta}{n+2+\xi+\eta} \quad \mbox{if} \quad n \quad \mbox{is} \quad \mbox{odd}  \right . , \quad
C_n= \left\{-\frac{n+1+2\xi}{n+1+\xi+\eta} \quad \mbox{if} \quad n \quad \mbox{is} \quad \mbox{even} \atop
-\lambda\frac{n}{n+1+\xi+\eta} \quad \mbox{if} \quad n \quad \mbox{is} \quad \mbox{odd}  \right ..
\ee

Set $\lambda=\frac{1+c}{1-c}$ and $g=-\frac{2}{1-c}$. Now,
introducing the new parameters \be A_n'=\frac{1}{g}A_n,\quad
C_n'=\frac{1}{g}C_n \ee and choosing the parametrization \be
\xi=\frac{\alpha-1}{2}, \quad \eta=\frac{\beta-1}{2}, \ee we can
easily see that the coefficients $A_n'$ and $C_n'$ correspond to
those of the big -1 Jacobi polynomials introduced in~\cite{VZ_big}
with a shifted spectral parameter. Recall that the big -1 Jacobi
polynomials are polynomials orthogonal on the union of two
intervals $[-1,-c]\cup [c,1]$ with respect to the weight function
\be w^{(-1)}(x) = \sign x(x+1)(x-c) (1-x^2)^{(\alpha-1)/2}(x^2
-c^2)^{(\beta-1)/2}. \lab{wbig} \ee

Notice that the case $\lambda=1$ corresponds to $c=0$. This special case corresponds to the little -1 Jacobi polynomials.
In this instance the two intervals of orthogonality connect to a form the
single interval $[-1,1]$.

The big -1 Jacobi polynomials $P_n^{(\alpha,\beta)}(x;c)$
corresponding to the weight function \re{wbig} are the
eigensolutions of \cite{VZ_big} \be L P_n^{(\alpha,\beta)}(x;c) =
\lambda_n P_n^{(\alpha,\beta)}(x;c), \lab{LPP} \ee where the
operator \be L^{(\alpha,\beta,c)}= g_0(x)(R-I) + g_1(x) \partial_x
R \lab{L_0} \ee with \be g_0(x)= \frac{(\alpha+\beta+1)x^2
+(c\alpha -\beta)x + c}{x^{2}}, \quad
g_1(x)=\frac{2(x-1)(x+c)}{x}, \lab{g_01} \ee $I$ is the identity
operator and $R$ the reflection operator $Rf(x) = f(-x)$.

The eigenvalues of $L^{(\alpha,\beta,c)}$ are \be \lambda_n =
\left\{ {2n, \quad n \quad \mbox{even}    \atop
-2(\alpha+\beta+n+1), \quad n \quad \mbox{odd}} \right . .
\lab{lam-1} \ee As shown in \cite{VZ_Bochner} the operator
\re{L_0} is the most general operator of the first order in
$\partial_x$ and involving $R$ that has orthogonal polynomials as
eigenfunctions.

We thus see that the orthogonal polynomials $Q_n(x;\lambda)$
obtained  by the SDG map from the Jacobi OPUC, satisfy a Dunkl type
eigenvalue equation. This equation can easily be obtained from
\re{LPP} by an affine transformation of the argument $x \to \mu x
+ \nu$ with appropriately chosen constants $\mu,\nu$.

Note that when $\xi=\eta=-1/2$ in \re{gen_ultra_a}, $a_n=0,
n=0,1,2,\dots$, i.e. this case is equivalent to the 1-periodic
case considered in Section~8. This corresponds to the case
$\alpha=\beta=0$ of the big -1 Jacobi polynomials. In particular,
the Chebyshev polynomials of the third kind \cite{MH}
(corresponding to $\lambda=1$,  or equivalently, $c=0$)
$$
V_n(x) = \frac{\cos(n(\tau+1/2))}{\cos(\tau/2)}, \quad x=\cos \tau
$$
coincide with the little -1 Jacobi polynomials \cite{VZ_little}
with $\alpha=\beta=0$ and hence they satisfy the Dunkl type
eigenvalue equation \be 2(x-1) V_n'(-x)+ V_n(-x) = (-1)^n
(2n+1)V_n(x). \lab{Ch_eig} \ee

Similarly, for the Chebyshev polynomials of the fourth kind
\cite{MH}
$$
W_n(x) = \frac{\sin(n(\tau+1/2))}{\sin(\tau/2)}, \quad x=\cos \tau
$$
one has the eigenvalue equation \be 2(x+1) W_n'(-x)+ W_n(-x) =
(-1)^n (2n+1)W_n(x). \lab{Ch4_eig} \ee In fact, equation
\re{Ch4_eig} is a simple consequence of \re{Ch_eig} because of the
property  $V_n(-x)=(-1)^n W_n(x)$  \cite{MH}.

\bigskip\bigskip
{\Large\bf Acknowledgments}

\bigskip

The research of LV is supported in part by a research grant from
the Natural Sciences and Engineering Research Council (NSERC) of
Canada.

\newpage

\end{document}